\documentclass[12pt,a4paper]{article}

\usepackage{amsmath}
\newcommand{\Ord}[1]{{\cal O}\big(#1\big)}

\newcommand{\half}{\mbox{$\frac12$}}

\newcommand{\rat}[2]{{\textstyle\frac{#1}{#2}}}
\def\theoremcolour{}
\newtheorem{theorem}{\theoremcolour Theorem}

\newtheorem{corollary}[theorem]{\theoremcolour Corollary}

\newenvironment{proof}%
{\par\textbf{\theoremcolour Proof:}}%
{\hfill{\theoremcolour$\spadesuit$}\par}
\newcommand{\hd}{h\partial}

\usepackage{url}

\IfFileExists{ajr.sty}{\usepackage[colour,hyperref]{ajr}}{}

\title{A corrected quadrature formula and applications} 

\author{Nenad
Ujevi\'c\thanks{Department of Mathematics, University of Split, Teslina
12/III, 21000 Split, \textsc{Croatia}.
\protect\url{mailto:ujevic@mapmf.pmfst.hr}} \and A. J.
Roberts\thanks{Dept Maths \& Computing, University of Southern
Queensland, Toowoomba, Queensland 4352, \textsc{Australia}.
\protect\url{mailto:aroberts@usq.edu.au}}}

\date{12 March 2003}

\begin{document}

    \maketitle

\begin{abstract}
A straightforward 3-point quadrature formula of closed type is
derived that improves on Simpson's rule.  Just using the additional
information of the integrand's derivative at the two endpoints we
show the error is sixth order in grid spacing.  Various error
bounds for the quadrature formula are obtained to quantify more
precisely the errors.  Applications in numerical integration are
given.  With these error bounds, which are generally better than
the usual Peano bounds, the composite formulas can be applied to
integrands with lower order derivatives.
\end{abstract}


\tableofcontents

\section{Introduction} 
\label{sec:intro}

In recent years some authors have considered so called perturbed
(corrected) quadrature rules.  For example, the corrected midpoint and
trapezoid quadrature rules are considered in~\cite{CD1}
and~\cite{CD2}.  As a specific introductory example, consider a
correction to the midpoint rule.  The classical midpoint rule has the
form
\begin{equation}
    \int_a^b f(x)\,dx = (b-a)f(\rat{a + b} 2 ) +R(f)\,, 
    \label{eq:1.1} 
\end{equation}
where $R(f)$ is the remainder term (error) of this rule.  A perturbed
(corrected) rule involves the derivative at the endpoints:
\begin{equation}
    \int_a^b f(x)\,dx = (b-a)f(\rat{a + b} 2 )
+ \frac{(b-a)^2}{ 24}\left[ f'(b)-f'(a)\right] + R_1(f)\,,  
\label{eq:1.2}
\end{equation}
where $R_1(f)$ is the remainder term (error) of this rule.  The
following properties have analogues in the work we report on Simpson's
rule:
\begin{enumerate}
	\item \label{en:1} the original rule~(\ref{eq:1.1}) is exact for
	polynomials of degree~1, while the perturbed rule is exact for
	polynomials of degree~3,

	\item \label{en:2} estimates of the errors are respectively
    \begin{eqnarray}
        |R(f)| &\leq& \frac{(b -a)^3}{ 24} M_2\,, \label{eq:1.3} \\
    |R_1(f)| &\leq& \frac{(b-a)^3}{18\sqrt3} M_2\,, \label{eq:1.4}
    \end{eqnarray}
	where $|f''(t)|\leq M_2$\,, $t \in [a, b]$ --- the
	bound~(\ref{eq:1.4}) is better than~(\ref{eq:1.3});

	\item a corresponding composite quadrature formula, for the
	corrected rule, has only one additional term, with respect to a
	composite formula for the original rule,

	\item \label{en:4} the corrected composite formula has a better
	estimation of error than the original composite formula (a
	consequence of~\ref{en:2}).
    
    Indeed another bound is 
\begin{equation*}
\left| R_{1}(f)\right| \leq \frac{7M_{4}}{5760}(b-a)^{5}\,,
\end{equation*}%
where $\left| f^{(4)}(t)\right| \leq M_{4}$ for $t\in [ a,b] $\,.

\end{enumerate}

The above properties are valid for similar corrected rules (for
example, for the corrected trapezoid rule)~\cite{CD2}.  However, we
cannot correct all quadrature rules such that \emph{all} the
properties~\ref{en:1}--\ref{en:4} hold.  In Section~\ref{sec:fd} we
show that the well-known Simpson's rule does have a simple endpoint
correction, but that the quadrature weights have to be modified as
well, see~(\ref{eq:mmsimp}).  We highlight some advantages of the
corrected rule over the Simpson's rule.  In Section~\ref{sec:main}
various error bounds of this rule are obtained.  These error bounds are
generally, but not always, better than the usual Peano error bounds.
In Section~\ref{sec:appl} applications in numerical integration are
given.  An illustrative example demonstrates that the modified
rule gives better results than Simpson's rule.

\section{Finite differences derive the modified Simpson's rule}
\label{sec:fd}

We modify Simpson's rule for integration.  First, we consider 
integration rules formed over just two consectutive subintervals, each 
of length~$h$, and derive endpoint modifications.  This is analogous 
to the improvement to~(\ref{eq:1.1}) achieved by the inclusion of 
endpoint derivative information in~(\ref{eq:1.2}).  Second, this 
modified Simpson's rule is straightforwardly summed to apply to an 
integration over many subintervals.  In later sections we rederive 
these formula with less restrictions on the integrand and with error 
bounds rather than just leading order estimates.

\begin{theorem}
	For $C^\infty[a,b]$ integrands~$f$,
		\begin{eqnarray}
		\int_a^bf(x)\,dx 
		&=&\frac{b-a}{30}\left[ 7f(a)+16f(\rat{a+b}2)+7f(b) \right]
        \nonumber\\&&{}
		-\frac{(b-a)^2}{60}\left[ f'(b)-f'(a) \right]
		+\bar R(f)\,,
		\label{eq:msimp}
	\end{eqnarray}
	where the error term is, to a leading order estimate,
		\begin{equation}
		\bar R(f)\approx \frac{(b-a)^6}{302400}\left[ f^v(b)-f^v(a) \right]\,.
		\label{eq:msimpe}
	\end{equation}
\end{theorem}

\paragraph{Example:} As a simple illustrative example, consider
$\int_{-1}^1 e^x\,dx=e-1/e=2.3504$\,.  Simpson's rule estimates the
integral as approximately $(e+4+1/e)/3=2.3621$\,, in error by
about~$0.01$\,.  However, the modified Simpson's rule~(\ref{eq:msimp})
estimates the integral as $(6e+16+8/e)/15=2.3502$\,, which has an error
about two orders of magnitude smaller.  This modification to Simpson's
rule can be very effective.

\begin{proof}
Consider integrating~$f$ over two consectutive intervals in a regular 
grid of points~$x_j$ with grid spacing~$h$.  Identify $a=x_{j-1}$, 
$b=x_{j+1}$ and hence $x_j=(a+b)/2$ is the midpoint.  Following 
\cite[p65]{npl61}, we write the analysis in terms of centred 
difference and mean operators, $\delta f_j=f_{j+1/2}+f_{j-1/2}$ and 
$\mu f_j=(f_{j+1/2}-f_{j-1/2})/2$ respectively, and the 
differentiation operator denoted by~$\partial$.  Then 
the intergral \cite[p69]{npl61}
\begin{displaymath}
\frac{1}{2h}\int_{x_{j-1}}^{x_{j+1}} f(x)\,dx	 = (\hd)^{-1}\mu\delta f_j  \,.
\end{displaymath}
To derive a three point integration rule with endpoint corrections 
such as~(\ref{eq:msimp}), the above right-hand side must be in the 
form $[1 +\alpha\delta^2 +\mu\delta \beta(\hd)]f_j$ for some 
constant~$\alpha$ and some function~$\beta(\hd)$: 
$[1+\alpha\delta^2]f_j$ symmetrically involves~$f_j$ and~$f_{j\pm1}$ 
alone; and $\mu\delta\beta(\hd)f_j$ only involves the derivatives 
of~$f$ at the endpoints~$x_{j\pm1}$\,.  Thus we rearrange the operator 
equation
\begin{eqnarray*}
	 &  & (\hd)^{-1}\mu\delta = 1 +\alpha\delta^2 +\mu\delta \beta(\hd) \\
	 & \Leftrightarrow & \beta(\hd)=\frac{1}{\hd} -\frac{1}{\mu\delta} 
	 -\alpha\frac{\delta}{\mu}  \\
	 & \Leftrightarrow & \beta(\hd)= \frac{1}{\hd} -\frac{1}{\sinh\hd} 
     -\alpha 2\tanh \half\hd\,,
\end{eqnarray*}
as $\delta=2\sinh\half\hd$ and $\mu=\cosh\half\hd$ \cite[p65]{npl61}.
Choosing this particular function~$\beta(\hd)$ would generate a rule in
the interior together with end point corrections that would give an
exact quadrature formula (for all~$\alpha$).  However, the infinite
derivatives required are not practical.

We chose~$\alpha$ to generate an accurate rule that only needs to know
function values and the end point derivatives.  The approach is to
expand this function~$\beta(\hd)$ in powers of small~$\hd$ to see
\begin{displaymath}
	\beta=(-\alpha+\rat16)\hd 
	+(\rat1{12}\alpha-\rat7{360})(\hd)^3 
	+(-\rat1{120}\alpha+\rat{31}{15120})(\hd)^5 +\cdots\,.
\end{displaymath}
Observe as an aside that choosing $\alpha=1/6$ eliminates the first
derivative term in~$\beta$ leading to the familiar Simpson's rule with
an error determined by the neglected parts of~$\beta$, namely the
end-point contributions $-\mu\delta(\hd)^3/180+\cdots$\,.\footnote{Thus
Simpson's rule has no simple end correction involving just the first
derivative of the integrand, of the type described in~(\ref{eq:1.2})
for the mid-point rule, instead the end correction would necessarily
involve the third derivative.} Instead we choose $\alpha=7/30$ to
eliminate the third derivative term in~$\beta$:
\begin{eqnarray}
	\int_{x_{j-1}}^{x_{j+1}}f(x)\,dx 
	& = & 2h(\hd)^{-1}\mu\delta f_j  \nonumber\\
	 & = &  2h\left[1+\rat7{30}\delta^2\right]f_j
	 +2h\mu\delta\left[-\rat1{15}\hd+\rat1{9450}\hd^5+\cdots\right]f_j\nonumber\\
	 & = & \rat h{15}\left[7f_{j-1}+16f_j+7f_{j+1}\right]
	 \nonumber\\&&{}
	 -\rat {h^2}{15}\left[f'_{j+1}-f'_{j-1}\right]
	 +\rat{h^6}{4725}\left[f^v_{j+1}-f^v_{j-1}\right]
	 +\cdots\,.
	 \label{eq:msimps}
\end{eqnarray}
Substitute $h=(b-a)/2$ to reproduce~(\ref{eq:msimp}) and its leading 
order error~(\ref{eq:msimpe}).
\end{proof}

\begin{corollary}
	Apply~(\ref{eq:msimps}) to~$n$ consecutive pairs of intervals 
	from say~$a=x_0$ to~$b=x_{2n}$ and sum to immediately deduce that for 
	$C^\infty[a,b]$ integrands~$f$
		\begin{equation}
		\int_a^bf(x)\,dx
		=\frac h{15}\sum_{j=1,\,j\text{ odd}}^{2n-1}
		\left[7f_{j-1}+16f_j+7f_{j+1}\right]
		-\frac {h^2}{15}\left[f'(b)-f'(a)\right]
		+\bar R(f)\,.
		\label{eq:mmsimp}
	\end{equation}
	where the error is to leading order
		\begin{equation}
		\bar R(f)\approx\frac{h^6}{4725}\left[f^v(b)-f^v(a)\right]\,.
		\label{eq:mmsimpe}
	\end{equation}
\end{corollary}

See that these simple modifications to Simpson's rule generate an 
integration method with error~$\Ord{h^6}$.

%
%
%
%
%

\section{Further error analysis}
\label{sec:main}

For the sake of simplicity, we first consider the error analysis for the
quadrature formula~(\ref{eq:msimp}) on the interval~$\left[ 0,1\right] $.
Then we easily transform obtained results to an arbitrary interval. The
formula~(\ref{eq:msimp}) on the interval~$\left[ 0,1\right] $ has the form
\begin{equation}
\int\limits_0^1f(x)\,dx=\frac{7f(0)+16f(\frac 12)+7f(1)}{30}
-\frac{f^{\prime}(1)-f^{\prime }(0)}{60}+R_k(f)\,,  \label{f1}
\end{equation}
for $k=2,3,4,5,6$\,. Using the Peano Kernel Theorem we find the following
Peano kernels
\begin{equation}
T_k(x)=\left\{ 
\begin{array}{ll}
(-1)^kP_k(x)\,,& x\in \left[ 0,\frac 12\right)\,; \\ 
(-1)^kQ_k(x)\,,& x\in \left[ \frac 12,1\right]\,;
\end{array}
\right.  \label{peano}
\end{equation}
for $k=2,3,4,5,6$\,, where
\begin{displaymath}
 \begin{array}{ll}
P_2(x)=\frac 12x^2-\frac 7{30}x+\frac 1{60}\,,
& Q_2(x)=\frac
12x^2-\frac{23}{30}x+\frac{17}{60}\,, \\ 
P_3(x)=\frac 1{3!}x(x-\frac 15)(x-\frac 12)\,,
& Q_3(x)=\frac
1{3!}(x-1)(x-\frac 12)(x-\frac 4{25})\,, \\ 
P_4(x)=\frac 1{4!}x^2(x-\frac 13)(x-\frac 35)\,,
& Q_4(x)=\frac
1{4!}(x-1)^2(x-\frac 23)(x-\frac 25)\,, \\ 
P_5(x)=\frac 1{5!}x^3(x-\frac 12)(x-\frac 23)\,,
& Q_5(x)=\frac
1{5!}(x-1)^3(x-\frac 12)(x-\frac 13)\,, \\ 
P_6(x)=\frac 1{6!}x^4(x^2-\frac 75x+\frac 12)\,,
& Q_6(x)=\frac
1{6!}(x-1)^4(x^2-\frac 35x+\frac 1{10})\,.
\end{array}
\end{displaymath}
We have
\begin{equation}
R_k(f)=\int_0^1T_k(x)f^{(k)}(x)\,dx\,,  \label{f2}
\end{equation}
for $k=2,3,4,5,6$\,. We also have 
\begin{equation}
\int_0^1T_k(x)dx=0\\,,\quad k=2,\ldots,5\,,  \label{a2}
\end{equation}
\begin{equation}
\int_0^1\left| T_k(x)\right| \,dx=C_k\\,\quad k=2,\ldots,6\,,  \label{a3}
\end{equation}
where 
\begin{eqnarray}&&
C_2=\frac{19\sqrt{19}}{10125}\,,\quad C_3=\frac{253}{360000}\,,\quad 
C_4=\frac 1{14580}\,,\quad C_5=\frac 1{115200}\,, \quad \label{a4}
\\&&
C_6=\frac 1{604800}\,,  \label{a5}
\end{eqnarray}
and
\begin{equation}
{\max_{x\in \left[ 0,1\right] } }\left| T_k(x)\right| =B_k\,,\quad
k=2,3,4,5\,,  \label{a6}
\end{equation}
where 
\begin{equation}
B_2=\frac 1{40}\,,\quad 
B_3=\frac 7{20250}+\frac{19\sqrt{19}}{81000}\,,\quad 
B_4=\frac 1{5760}\,,\quad 
B_5=\frac 1{58320}\,.  \label{a7}
\end{equation}

\begin{theorem}
\label{T2} Let $f\in C^k(0,1)$ and let $\gamma _k$, $\Gamma _k$ be real
numbers such that $\gamma _k\leq f^{(k)}(x)\leq \Gamma _k$\,, $x\in
\left[ 0,1\right] $, $k=2,3,4,5$\,.  Let $S_k=f^{(k)}(1)-f^{(k)}(0)$\,,
$k=1,2,3,4$\,.  Then we have
\begin{eqnarray}&&
\left| R_k(f)\right| \leq \frac{\Gamma _k-\gamma _k}2C_k\,,\quad 
k=2,\ldots,5\,,
\label{a12}
\\&&
\left| R_k(f)\right| \leq (S_{k-1}-\gamma _k)B_k\,,\quad k=2,\ldots,5\,,
\label{a13}
\\&&
\left| R_k(f)\right| \leq (\Gamma _k-S_{k-1})B_k\,,\quad k=2,\ldots,5\,,
\label{a14}
\end{eqnarray}
where $R_k(f)$ are defined by~(\ref{f2}), $C_k$ are defined by~(\ref{a4})
and $B_k$ are defined by~(\ref{a7}), for $k=2,3,4,5$\,.
\end{theorem}

\begin{proof}
Let $C$ be an arbitrary constant. Then we have 
\begin{equation}
R_k(f)=\int_0^1T_k(x)\left[ f^{(k)}(x)-C\right] \,dx
=\int_0^1T_k(x)f^{(k)}(x)\,dx\,,  \label{a15}
\end{equation}
for $k=2,3,4,5$\,, since~(\ref{a2}) holds.
If we now choose $C=({\Gamma _k+\gamma _k})/2$\,, then we get 
\begin{eqnarray}
\left| R_k(f)\right| &=&\left| \int_0^1T_k(x)\left[ f^{(k)}(x)-\frac{
\Gamma _k+\gamma _k}2\right] \,dx\right|  \label{a16} \\
&\leq& {\sup_{x\in \left[ 0,1\right] } }\left| f^{(k)}(x)-\frac{
\Gamma _k+\gamma _k}2\right| \int_0^1\left| T_k(x)\right| \,dx  \notag
\\
&\leq& \frac{\Gamma _k-\gamma _k}2\int_0^1\left| T_k(x)\right| \,dx
\nonumber\\&=&
\frac{\Gamma _k-\gamma _k}2C_k\,,  \notag
\end{eqnarray}
for $k=2,3,4,5$\,.
If we choose $C=\gamma _k$\,, then we have 
\begin{eqnarray*}
\left| R_k(f)\right| &=&\left| \int_0^1T_k(x)\left[ f^{(k)}(x)-\gamma _k
\right] \,dx\right| \\
&\leq& {\max_{x\in \left[ 0,1\right] } }\left| T_k(x)\right|
\int_0^1\left| f^{(k)}(x)-\gamma _k\right| \,dx \\
&=&B_k\int_0^1\left[ f^{(k)}(x)-\gamma _k\right] \,dx
\\&=&(S_{k-1}-\gamma_k)B_k\,,
\end{eqnarray*}
for $k=2,3,4,5$\,.

In a similar way we can prove that~(\ref{a14}) holds.
\end{proof}

The estimations~(\ref{a12}) are Peano-like bounds and they are
generally (but not always) better than the usual Peano bounds.  Namely,
we know that~$T_k(x)$, $k=2,\ldots,6$ are Peano kernels.  The usual
Peano error bounds are
\begin{equation*}
\left| R_k(f)\right| \leq \left\| f^{(k)}\right\| _\infty
\int_0^1\left| T_k(x)\right| \,dx=C_k\left\| f^{(k)}\right\| _\infty 
\,,\quad k=2,\ldots,6\,,
\end{equation*}
where $\left\| f^{(k)}\right\| _\infty ={ \sup_{x\in \left[ 0,1\right]
} }\left| f^{(k)}(x)\right| $\,.  If we choose $\gamma _k={\inf_{x\in
\left[ 0,1\right] } }f^{(k)}(x)$ and $\Gamma _k={\sup_{x\in \left[
0,1\right] } }f^{(k)}(x)$\,, then $\frac{\Gamma _k-\gamma _k}2\leq
\left\|f^{(k)} \right\| _\infty $\,.  Thus, in this case, the error
bounds given by~(\ref{a12}) are better than the Peano error bounds.  In
fact, they are equal if and only if $\Gamma _k=-\gamma _k$\,.  This
case ($\Gamma _k=-\gamma _k$) is very rare in practice.

Theoretically, we can derive better error bounds.  Let us say something
about the last assertion.
We can verify that 
\begin{equation*}
\int_0^1T_k(x)p_j(x)\,dx=0\,,\quad j=5-k\,,\quad k=2,3,4\,,
\end{equation*}
where $p_j(x)$ is any polynomial of degree~$\leq j$\,. Thus, 
\begin{equation*}
\int_0^1T_k(x)\left[ f^{(k)}(x)-p_j(x)\right] \,dx=\int
_0^1T_k(x)f^{(k)}(x)\,dx
\end{equation*}
such that 
\begin{equation*}
\left| R_k(f)\right| \leq \left\| f^{(k)}-p_j \right\| _\infty C_k
\,,\quad 
k=2,3,4\,,\quad j=5-k\,.
\end{equation*}
The above estimations are theoretically better than the
corresponding estimations in~(\ref{a12}).

We now give the above obtained results for an arbitrary
interval~$\left[ a,b \right] $.  The mapping $x=({t-a})/({b-a})$ is a
bijection from~$\left[ a,b \right] $ onto~$\left[ 0,1\right] $.  If we
use this bijection then we find that the polynomials~$P_k$ and~$Q_k$ on
the interval~$\left[ a,b\right] $ have the forms: $\tilde P_0(t)=1$\,,
$\tilde Q_0(t)=1$\,, $\tilde P_1(t)=t-({ 23a+7b})/{30}$\,, $\tilde
Q_1(t) =t-({7a+23b})/{30}$\,, etc.  The polynomials can be also
obtained by simple integration.  For example, $\tilde P_2$ can be
obtained by integrating~$\tilde P_1$ and determining an additional
constant such that~(\ref{a2}) holds.  (They are additionally
normalized.)  We define the functions
\begin{equation}
\tilde T_k(t)=\left\{ 
\begin{array}{ll}
(-1)^k\tilde P_k(t)\,,& t\in \left[ a,\frac{a+b}2\right)\,; \\ 
(-1)^k\tilde Q_k(t)\,,& t\in \left[ \frac{a+b}2,b\right]\,;
\end{array}
\right.  \label{a1a}
\end{equation}
for $k=2,\ldots,6$\,.  Here we choose $a=x_{i-1}$\,, $b=x_{i+1}$ and
use the notations introduced in Section~\ref{sec:fd}.  Using these
functions we derive the following results.

\begin{corollary}
 \label{C1}Let $f\in C^{k}(x_{j-1},x_{j+1})$\,, $k=2,3,\ldots,6$\,. Then we
have 
\begin{equation}
\int_{x_{j-1}}^{x_{j+1}}f(x)\,dx=\frac{7f_{j-1}+16f_{j}+7f_{j+1}}{15}h-
\frac{f_{j+1}^{\prime }-f_{j-1}^{\prime}}{15}h^{2}
+\tilde{R}_{k}(f)\,,
\label{a8a}
\end{equation}
where 
\begin{equation}
\tilde{R}_{k}(f)=\int_{x_{j-1}}^{x_{j+1}}\tilde{T}_{k}(x)f^{(k)}(x)\,dx
\\,,  \label{a9a}
\end{equation}
for $k=2,\ldots,6$ and $\tilde{T}_{k}(x)$ defined by~(\ref{a1a}).
\end{corollary}

\begin{corollary}
\label{C2} Let $f\in C^{k}(x_{j-1},x_{j+1})$ and let $\gamma _{k}$,
$\Gamma _{k}$ be real numbers such that $\gamma _{k}\leq f^{(k)}(x)\leq
\Gamma _{k}$\,, $x\in \left[ x_{j-1},x_{j+1}\right] $\,, $k=2,3,4,5$\,.
Let 
\begin{displaymath}
    S_{k}=\frac{
f^{(k)}(x_{j+1})-f^{(k)}(x_{j-1})} {x_{j+1}-x_{j-1}}\,,\quad
k=1,2,3,4\,.
\end{displaymath}
Then we have
\begin{eqnarray}&&
\left| \tilde{R}_{k}(f)\right| \leq \frac{\Gamma _{k}-\gamma _{k}}{2}
D_{k}h^{k+1}\,,  \label{a12a}
\\&&
\left| \tilde{R}_{k}(f)\right| \leq (S_{k-1}-\gamma _{k})E_{k}h^{k+1}\,,
\label{a13a}
\\&&
\left| \tilde{R}_{k}(f)\right| \leq (\Gamma _{k}-S_{k-1})E_{k}h^{k+1}\,,
\label{a14a}
\end{eqnarray}
where~$\tilde{R}_{k}(f)$ are defined by~(\ref{a1a}),
$D_{k}=2^{k+1}C_{k}$ ($ C_{k}$ are defined by~(\ref{a4})) and
$E_{k}=2^{k+1}B_{k}$ ($B_{k}$ are defined by~(\ref{a7})),
for~$k=2,3,4,5$\,.
\end{corollary}

\begin{corollary}
\label{C3}Let $f\in C^{6}(x_{j-1},x_{j+1})$\,. Then we have 
\begin{equation}
\left| \tilde{R}_{6}(f)\right| \leq D_{6}\left\| f^{(6)}\right\| _{\infty
}h^{7} \,, \label{a15a}
\end{equation}
where $D_{6}=2^{7}C_{6}$ ($C_{6}$ is defined by~(\ref{a5})) and $\left\|
g\right\| _{\infty }={\sup_{x\in \left[ x_{j-1},x_{j+1}\right] } }
\left| g(x)\right| $\,.
\end{corollary}

\section{Applications in numerical integration}
\label{sec:appl}

We define the partition $\pi =\left\{ a=x_0<x_1<\cdots
<x_{2n}=b\right\}$ of the interval~$\left[ a,b\right] $ such that
$x_{i+1}=x_i+h$\,, $i=0,1,\ldots,2n-1$\,, $h=(b-a)/(2n)$\,.  We also
define the functions
\begin{equation*}
\bar T_{ki}(x)=\left\{ 
\begin{array}{ll}
(-1)^k\bar P_{ki}(x)\,,& x\in \left[ x_{i-1},x_i\right] \,; \\ 
(-1)^k\bar Q_{ki}(x)\,,& x\in \left( x_i,x_{i+1}\right] \,;
\end{array}
\right.
\end{equation*}
for $k=2,\ldots,6$\,, $i=1,3,\ldots,2n-1$\,, which correspond to the
functions defined by~(\ref{a1a}) (for~$a=x_{i-1}$, $b=x_{i+1}$).

\begin{theorem}
 \label{T31} Under the assumptions of Corollary~\ref{C1} suppose that~$
\pi $ and~$\bar{T}_{k}$ are given as above. Then we have 
\begin{eqnarray}
\int_{a}^{b}f(x)\,dx &=&\frac{h}{15}\sum_{j=1,j \text{
odd}}^{2n-1}\left[ 7f_{j-1}+16f_{j}+7f_{j+1}\right] \label{g1} \\
&&{} -\frac{h^{2}}{15}\left[ f^{\prime }(b)-f^{\prime }(a)\right] +\bar{R}
_{k}(f)\,,  \notag
\end{eqnarray}
where 
\begin{equation}
\bar{R}_{k}(f)=\sum_{j=1,\,j\text{ odd}}^{2n-1}\int
_{x_{j-1}}^{x_{j+1}}\bar{T}_{ki}(x)f^{(k)}(x)\,dx\,,  \label{g2}
\end{equation}
for~$k=2,\ldots,6$\,.
\end{theorem}

\begin{proof}
Sum~(\ref{a8a}) over odd~$j$ from~$1$ to~$2n-1$ to get~(\ref{g1}--\ref{g2}),
since
\begin{equation*}
\sum_{j=1,\,j\text{ odd}}^{2n-1}\left[ f_{j+1}^{\prime
}-f_{j-1}^{\prime }\right] =f^{\prime }(b)-f^{\prime }(a)\,.
\end{equation*}
\end{proof}

\begin{theorem}
 \label{T32}Under the assumptions of Theorem~\ref{T31} and
 Corollary~\ref{C2} we have
\begin{eqnarray*}&&
\left| \bar R_k(f)\right| \leq \frac{\Gamma _k-\gamma _k}4D_kh^k(b-a)\,,
\\&&
\left| \bar R_k(f)\right| \leq \frac 12(S_{k-1}-\gamma _k)E_kh^k(b-a)\,,
\\&&
\left| \bar R_k(f)\right| \leq \frac 12(\Gamma _k-S_{k-1})E_kh^k(b-a)\,,
\end{eqnarray*}
for~$k=2,\ldots,5$\,.
\end{theorem}

\begin{proof}
The proof follows immediately from Theorem~\ref{T31} and
Corollary~\ref{C2}.
\end{proof}

\begin{theorem}
 Under the assumptions of Theorem~\ref{T31} and Corollary~\ref{C3} we
have 
\begin{equation*}
\left| \bar R_6(f)\right| \leq \frac 12D_6h^6(b-a)\left\| f^{(6)}\right\|
_\infty \,.
\end{equation*}
\end{theorem}

\begin{proof}
The proof follows immediately from Theorem~\ref{T31} and Corollary~\ref{C3}.
\end{proof}

Finally, let us compare the rule obtained in Theorem~\ref{T31} with the
standard composite Simpson's rule 
\begin{equation}
\int_{a}^{b}f(x)\,dx=\frac{h}{3}\sum_{j=1,\,j\text{ odd}
}^{2n-1}\left[ f_{j-1}+4f_{j}+f_{j+1}\right] +R_{S}(f)\,.  \label{sr}
\end{equation}
The terms 
\begin{equation*}
\frac{h}{3}\sum_{j=1,\,j\text{ odd}}^{2n-1}\left[
f_{j-1}+4f_{j}+f_{j+1}\right]\,,
\end{equation*}
and 
\begin{equation*}
\frac{h}{15}\sum_{j=1,\,j\text{ odd}}^{2n-1}\left[
7f_{j-1}+16f_{j}+7f_{j+1}\right]\,,
\end{equation*}
require a same amount of calculations. The rule~(\ref{g1}) has only one
additional term with respect to the rule~(\ref{sr}), namely 
\begin{equation*}
-\frac{h^{2}}{15}\left[ f^{\prime }(b)-f^{\prime }(a)\right] \,.
\end{equation*}
Hence, the amount of calculations is approximately the same for both
formulae.  (Recall that function evaluations are generally considered
the computationally most expensive part of quadrature algorithms.)  On
the other hand, the rule~(\ref{g1}) is exact for polynomials of
degree~$ \leq 5$, while the rule~(\ref{sr}) is exact for polynomials of
degree~$\leq 3 $\,.  Furthermore, from Theorem~\ref{T32} we have
\begin{equation}
\left| \bar R_4(f)\right| \leq 2\frac{\Gamma _4-\gamma _4}{3645}h^4(b-a)\,,
\label{e1}
\end{equation}
while the standard estimation for the Simpson's rule is 
\begin{equation}
\left| R_S(f)\right| \leq \frac{\left\| f^{(4)}\right\| _\infty
}{180}
h^4(b-a)\,.  \label{e2}
\end{equation}
Since, (\ref{e1}) is better than~(\ref{e2}), the rule~(\ref{g1}) has
better approximation properties than Simpson's rule.  Thus, we 
expect that it will give better results in practice (in most cases).

\begin{figure}
    \centering  
    \begin{tabular}{cc}
        \rotatebox{90}{\hspace{10ex}integration error} & \includegraphics{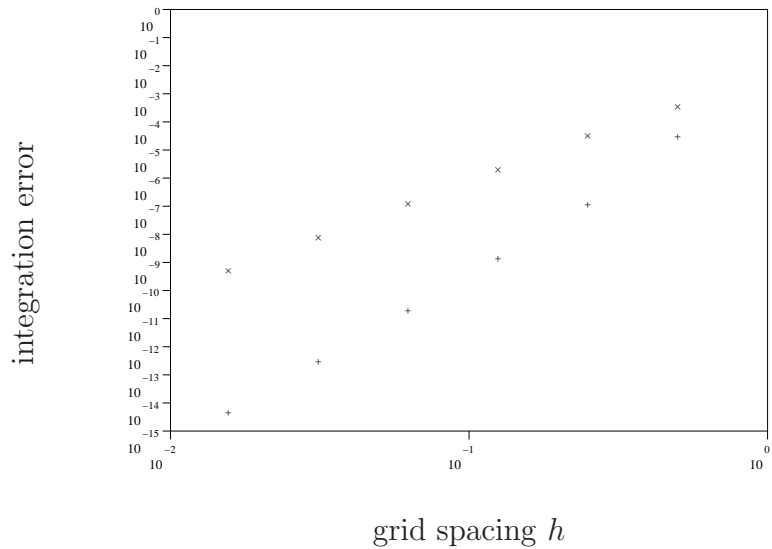}  \\
         & grid spacing $h$  
    \end{tabular}
	\caption{log-log plot of the errors of our integration
	rule~(\ref{eq:mmsimp}), $+$'s, as a function showing the
	$\Ord{h^6}$ rate of convergence to the erf
	integral~(\ref{eq:myerf}), compared to the $\Ord{h^4}$~convergence
	of the normal Simpson's rule, $\times$'s.}
    \label{fig:erf}
\end{figure}

\paragraph{Example:} here we show errors in estimating
\begin{equation}
    I=\frac{\sqrt\pi}2\mbox{erf}(1)
    =\int_0^1\exp(-x^2)\,dx\,.
    \label{eq:myerf}
\end{equation}
Using just two subintervals, $h=1/2$, our formula~(\ref{eq:mmsimp})
computes $I\approx 0.746795$\,, whereas with four subintervals, that
is
$h=1/4$, (\ref{eq:mmsimp}) gives $I\approx 0.746824$ which is correct
to six decimal places.  Figure~\ref{fig:erf} shows our rule converges
quickly with decreasing grid spacing~$h$, and is essentially exact to
double precision with just 64~subintervals.


\begin{thebibliography}{99}

\bibitem{CD1} P. Cerone and S. S. Dragomir, Midpoint-type Rules from an
Inequalities Point of View, Handbook of Analytic-Computational Methods in
Applied Mathematics, Editor: G. Anastassiou, CRC Press, New York, (2000),
135--200.

\bibitem{CD2} P. Cerone and S. S. Dragomir, Trapezoidal-type Rules from an
Inequalities Point of View, Handbook of Analytic-Computational Methods in
Applied Mathematics, Editor: G. Anastassiou, CRC Press, New York, (2000),
65--134.

%
%
%
%
%

\bibitem{npl61} National Physical Laboratory, Modern Computing Methods,
volume~16 of Notes on Applied Science, Her Majesty's Stationary Office, 1961.


\end{thebibliography}
\end{document}